\documentstyle[12pt]{article}
\font\teneufm=eufm10 \font\seveneufm=eufm7 \font\fiveeufm=eufm5
\newfam\frakturfam
\textfont\frakturfam=\teneufm \scriptfont\frakturfam=\seveneufm
\scriptscriptfont\frakturfam=\fiveeufm

\newtheorem{pr}{Proposition}

\newtheorem{lm}{Lemma}
\newtheorem{th}{Theorem}
\newtheorem{co}{Corollary}

\newtheorem{prob}{Problem}
\def\bee{\begin{eqnarray}}
\def\bes{\begin{eqnarray*}}
\def\eee{\end{eqnarray}}
\def\ees{\end{eqnarray*}}

\def\a{\alpha}

\def\Proof{{\sl Proof.}\ }

\title{Automorphisms and derivations of free Poisson algebras in two variables}

\begin{document}
\date{}
\maketitle

\begin{center}
{\textbf{Leonid Makar-Limanov} \footnote{Supported
by an NSA grant and by grant FAPESP, processo  06/59114-1.}}\\
Department of Mathematics \& Computer Science,\\
Bar-Ilan University, 52900 Ramat-Gan, Israel and\\
Department of Mathematics, Wayne State University, \\
Detroit, MI 48202, USA\\
e-mail: {\em lml@math.wayne.edu}\\

and\\

{\bf Umut Turusbekova, \ \ Ualbai Umirbaev}\\
Department of Mathematics, Eurasian National University\\
 Astana, 010008, Kazakhstan \\
e-mail: {\em umut.math@mail.ru, \ \ umirbaev@yahoo.com}

\end{center}

\begin{abstract}
Let $P$ be a free Poisson algebra in two variables over a field
of characteristic zero. We prove that the automorphisms of $P$ are
tame and that the locally nilpotent derivations of $P$ are
triangulable.
\end{abstract}

\noindent {\bf Mathematics Subject Classification (2000):} Primary
17B63, 14H37; Secondary 17B40, 17A36, 16W20.

\noindent

{\bf Key words:} Poisson algebras, automorphisms, derivations.

\section{Introduction}

\hspace*{\parindent}

 It is well known \cite{Czer,Jung,Kulk,Makar} that the automorphisms
of polynomial algebras and free associative algebras in two
variables are tame. It was recently proved \cite{Umi25,Um19} that
polynomial algebras and free associative algebras in three
variables in the case of characteristic zero have wild
automorphisms. P. Cohn  \cite{Cohn2} proved that the automorphisms
of a free Lie algebra with a finite set of generators are tame.

There are many other results, some of them quite deep, known about
the structure of polynomial algebras, free associative algebras,
and free Lie algebras. Though free Poisson algebras are very
closely connected with these algebras, only few results are known
about them up to now. Say, one of the fundamental results about
free associative algebras is the Bergman Centralizer Theorem (see
\cite{Berg}) which says that the centralizer of any nonconstant
element is a polynomial algebra on a single variable. An analogue of
this theorem for free Poisson algebras in the case of
characteristic zero was proved in \cite{MakarU2}.

The question on the tameness of automorphisms of free Poisson
algebras in two variables was open and was formulated in
\cite[Problem 5]{MakarU2}.  Note that the Nagata automorphism
\cite{Nagata,Umi25} gives an example of a wild automorphism of a
free Poisson algebra in three variables.

In \cite{Rentschler} R.\,Rentschler proved that the locally
nilpotent derivations of polynomial algebras in two variables over
a field of characteristic $0$ are triangulable. Using this result
he gave a new proof of Jung's Theorem \cite{Jung} on the tameness
of automorphisms of these algebras.

In this paper we study automorphisms and locally nilpotent
derivations of free Poisson algebras over a field of
characteristic zero. In Section 2 we introduce several gradings of
free Poisson algebras and describe some properties of homogeneous
derivations of these algebras. In Section 3 we prove that the
locally nilpotent derivations of two generated free Poisson
algebras are triangulable and the automorphisms of these algebras
are tame. These results are analogues of Rentschler's Theorem
\cite{Rentschler} and Jung's Theorem \cite{Jung}, respectively.

\section{Homogeneous derivations}

\hspace*{\parindent}

A vector space $B$ over a field $k$ endowed with two bilinear
operations $x\cdot y$ (a multiplication) and $\{x,y\}$ (a Poisson
bracket) is called {\em a Poisson algebra} if $B$ is a commutative
associative algebra under $x\cdot y$, $B$ is a Lie algebra under
$\{x,y\}$, and $B$ satisfies the following identity (the Leibniz
identity): \bes \{x, y\cdot z\}=\{x,y\}\cdot z + y\cdot \{x,z\}.
\ees Of course, the Leibniz identity just says that for every
$x\in B$ the map \bes ad_x : B\longrightarrow B, \ \ (y\mapsto
\{x,y\}), \ees is a derivation of $B$ as an associative algebra.

The map $ad_x$ also satisfies another similar identity: \bes
ad_x\{y,z\} =\{ad_x(y),z\} +  \{y,ad_x(z)\}. \ees It is just the
Jacobi identity for $B$ as a Lie algebra.

Let us call a linear homomorphism $D$ of $B$ to $B$ a derivation
of $B$ as a Poisson algebra if it satisfies both the Leibniz and
Jacobi identities. In other words, $D$ is simultaneously a
derivation of $B$ as an associative algebra and as a Lie algebra.

There are two important classes of Poisson algebras.

1) Symplectic algebras $S_n$. For each $n$ algebra $S_n$ is a polynomial algebra
$k[x_1,y_1, \ldots,x_n,y_n]$ endowed with the Poisson bracket
defined by \bes \{x_i,y_j\}=\delta_{ij}, \ \ \{x_i,x_j\}=0, \ \
\{y_i,y_j\}=0, \ees where $\delta_{ij}$ is the Kronecker symbol
and $1\leq i,j\leq n$.

2) Algebras of Lie type. Let $g$ be a Lie algebra with a linear
basis $e_1,e_2,\ldots,e_k,\ldots$. The symmetric algebra $S(g)$ of
$g$ (i. e. the usual polynomial algebra $k[e_1,
e_2,\ldots,e_k,\ldots]$) endowed with the Poisson bracket defined
by  \bes \{e_i,e_j\}=[e_i,e_j] \ees for all $i,j$, where $[x,y]$
is the multiplication of the Lie algebra $g$ is the Poisson algebra of type $g$.

From now on let $g$ be a free Lie algebra with free (Lie)
generators $x_1,x_2,\ldots,x_n$. It is well known (see, for
example \cite{Shest}) that in this case $S(g)$ is a free Poisson
algebra on the same set of generators. We denote this algebra
by $P=P\langle x_1,x_2,\ldots,x_n\rangle$. 

By $\deg$ we denote the
standard degree function of the homogeneous algebra $P$, i.e.
$\deg(x_i)=1$, where $1\leq i \leq n$. Note that
\begin{eqnarray*}
\deg\,\{f,g\}= \deg\,f+\deg\,g
\end{eqnarray*}
if $f$ and $g$ are homogeneous and $\{f,g\}\neq 0$. By
$\deg_{x_i}$ we denote the degree function on $P$ with respect to
$x_i$. We have $\deg_{x_i}(x_j)=\delta_{ij}$, where $1\leq i,j
\leq n$. The homogeneous elements of $P$ with respect to
$\deg_{x_i}$ can be defined in the ordinary way.

If $f$ is homogeneous with respect to each $\deg_{x_i}$, where
$1\leq i \leq n$, then $f$ is called multihomogeneous. For every
multihomogeneous element $f \in P$ we put \bes
mdeg(f)=(m_1,m_2,\ldots,m_n), \ees where $\deg_{x_i}f=m_i$ for all
$i$ and $1\leq i \leq n$.

Let us choose a multihomogeneous linear basis \bes
x_1,x_2,\ldots,x_n,\,[x_1,x_2],\ldots,[x_1,x_n],\ldots,[x_{n-1},x_n],\,[[x_1,x_2],x_3],\ldots
\ees of the free Lie algebra $g$ and denote the elements of this
basis by \bee\label{f1} e_1, e_2, \ldots, e_m, \ldots. \eee Note
that \bes mdeg \{e_i, e_j\} = mdeg(e_i) + mdeg(e_j) \ees if $i\neq
j$. So if $i < j$ then $\{e_i, e_j\}$ is a linear combination of
$e_m$ where all $m > j$.

The algebra $P=P\langle x_1,x_2,\ldots,x_n\rangle$ coincides with
the polynomial algebra on the elements (\ref{f1}). Consequently,
the words  \bee\label{f2} u=e_{i_1}e_{i_2}\ldots e_{i_k}, \ \
i_1\leq i_2\leq \ldots \leq i_k \eee form a linear basis of $P$.
The basis (\ref{f2}) is multihomogeneous since so is (\ref{f1}).

Consider the Lie algebra $Der(P)$ of all derivations of the
Poisson algebra $P$. For every system of elements
$f_1,f_2,\ldots,f_n$ of $P$ denote by \bee\label{f3}
D=\sum_{i=1}^{n} f_i \frac{\partial}{\partial x_i} \eee a unique
derivation  of $P$ such that $D(x_i)=f_i$ where $1\leq
i\leq n$. Then the derivations  \bee\label{f4}
v=u\frac{\partial}{\partial x_i}, \eee where $1\leq i\leq n$ and $u$ is 
an element of (\ref{f2}), constitute a linear basis
of $Der(P)$. For every element $v$ of the form (\ref{f4}) we put
\bes mdeg(v)=mdeg(u)-\epsilon_i, \ees where $\epsilon_i\in Z^n$ is
the standard basis vector with 1 in the $i$th position and with
zeroes everywhere else. Now one can define the multihomogeneous
derivations of the algebra $P$ and every element of $Der(P)$ can
be uniquely represented as the sum of multihomogeneous derivations
of different multidegrees.

To each nonzero vector $w\in Z^n$ we associate the so called
$w$-degree (or weight degree) function $wdeg$ on $P$ and $Der(P)$.
Put \bes wdeg(u)=<mdeg(u),w>, \ \ wdeg(v)=<mdeg(v),w>, \ees where
$u$ and $v$ are elements of the form (\ref{f2}) and (\ref{f4})
respectively, and  $<\,,\,>$ is the standard inner product in
$R^n$. Let $P_m$ and $Der_mP$ be the subsets of all
$w$-homogeneous elements of degree $m$ of $P$ and $Der(P)$,
respectively. It is clear that the decompositions \bes
P=\oplus_{m\in Z}P_m, \ \ Der(P)=\oplus_{m\in Z}Der_mP \ees are
gradings of the corresponding algebras. Moreover, for every
element $d\in Der_mP$ we have \bes d(P_k)\subseteq P_{m+k}. \ees

There is another natural degree function on $P$, just the total
degree  on $P$ as a polynomial ring, where the degree is one for
all elements of the homogeneous basis (\ref{f1}). Denote it by $pdeg$
and observe that
\bes
pdeg[a, b] = pdeg a + pdeg b - 1
\ees
 for any
$p$-homogeneous $a,\, b \in P$ if $[a, b] \neq 0$.

If $v$ is an element of the form (\ref{f4}) then we put \bes pdeg
v = pdeg u -1. \ees Let $P_m^*$ and $Der_m^*P$ be the subsets of
all $p$-homogeneous elements of degree $m$ of $P$ and $Der(P)$,
respectively. It is again clear that the decompositions \bes
P=\oplus_{m\in Z}P_m^*, \ \ Der(P)=\oplus_{m\in Z}Der_m^*P \ees
are gradings of the corresponding algebras and that for every
element $d\in Der_m^*P$ we have \bes d(P_k^*)\subseteq P_{m+k}^*.
\ees

Recall that a derivation $D$ of an algebra $R$ is called locally
nilpotent if for every $a\in R$ there exists a natural number
$m=m(a)$ such that $D^m(a)=0$. The statement of the next
proposition is well known (see, for example \cite[Proposition
5.1.15]{Essen}).

\begin{pr}\label{pr1}
Let $R=\oplus_{m\in Z}R_m$ be a graded algebra and suppose $D$ be a locally nilpotent derivation of $R$ such that
\bes
D=D_p+D_{p+1}+\ldots+D_q, \ \ D_i(R_m)\subseteq R_{i+m}, \ \ p\leq i\leq q, \ \ D_q\neq 0.
\ees
 Then $D_q$ is locally nilpotent.
\end{pr}
\Proof If
\bes
f=f_r+f_{r+1}+\ldots+f_s \in R,
\ees
where $f_i\in R_i$, $r\leq i\leq s$, and $f_s\neq 0$, then we put $\widehat{f}=f_s$.

Let $a\in R_m$ and assume that $D_q^i(a)\neq 0$ for any $i$. It
can be easily proved by induction on $i$ that \bes
\widehat{D^i(a)}=D_q^i(a). \ees Consequently, $D^i(a)\neq 0$ for
any $i$ and this gives a contradiction. $\Box$

Let $f$ be an arbitrary element of $P$ and $D$ be an arbitrary
derivation of $P$ of the form (\ref{f3}). We put \bes f
D=\sum_{i=1}^{n} (f f_i) \frac{\partial}{\partial x_i}. \ees Put
also \bes S(f)=\{e_{i_1}, e_{i_2},\ldots,e_{i_k}\} \ees if $f\in
k[S(f)]$ and $f\notin k[S(f_i)\setminus \{e_{i_j}\}]$, where
$1\leq j\leq k$. For $D$ we put \bes S(D)=S(f_1)\cup S(f_2)\cup
\ldots \cup S(f_n). \ees

If $x=e_i$ then we denote by $pdeg_x$ the polynomial degree function  with respect to $x$ on $P$.  Elements
$f\in P$ and $D\in Der\,P$ can be uniquely written as
\bes
f=f_0+xf_1+\ldots+x^mf_m, \ \ x\notin S(f_i),  \ \ 0\leq i \leq m,
\ees
and
\bes
D=D_0+xD_1+\ldots+x^mD_m, \ \ x\notin S(D_i), \ \ 0\leq i\leq m,
\ees
respectively. If $f_m\neq 0$ then $pdeg_x(f)=m$ and we put $l_x(f)=f_m$.
Put also $pdeg_xD=m$ and $l_x(D)=D_m$ if $D_m\neq 0$.

Put $e_i < e_j$ if $i < j$.

\begin{pr}\label{pr2}
Let $D$ be a derivation of $P$ and $x$ be the minimal element of
$S(D)$. Then \bes pdeg_xD(f)\leq pdeg_xD+pdeg_xf. \ees This
inequality becomes an equality iff $l_x(D)(l_x(f))\neq 0$ and in
this case \bes l_x(D(f))=l_x(D)(l_x(f)). \ees
\end{pr}
\Proof  Without loss of generality we may assume that $f$ is an
element of the basis (\ref{f2}) and $D$ is an element of the basis (\ref{f4}). 

If $f=uv$ then $D(f)=D(u)v+uD(v)$. So if the Proposition is true for $u$ and $v$ 
it is also true for $f$. Because of that we can assume that the polynomial degree 
of $f$ is one. Let us prove that in this case $pdeg_xD(f) \leq pdeg_xD$.

If $f=\{u,v\}$ then $D(f)=\{D(u),v\} + \{u,D(v)\}$. Denote by $L(x)$ the set of all 
elements $e_i$ such that  $e_i > x$.
If, say $D(u) = x^d u_1$ where $S(u_1) \subset L(x)$ then $\{D(u),v\} = \{x^d u_1,v\} = x^d \{u_1,v\} + dx^{d-1} u_1\{x,v\}$. As we remarked if $i < j$ then $\{e_i, e_j\}$ is a linear combination of $e_m$ where all $m > j$. So both $S(\{ u_1,v\})$ and $S(\{x,v\})$ are subsets of $L(x)$ and we can conclude that $pdeg_xD(f) \leq pdeg_xD$ if it is true for $u$ and $v$. It remains to check that $pdeg_xD(f) \leq pdeg_xD$ for $f$ with $deg(f) = 1$. Since we can assume that $D = x^du\frac{\partial}{\partial x_i}$ where $u \in L(x)$ and $1 \leq i \leq n$ we have $D(x_j) = 0$ when $j \neq i$ and $D(x_i) = x^du$.

So we proved that $pdeg_xD(f) \leq pdeg_xD + pdeg_xf$. To prove that $l_x(D(f))=l_x(D)(l_x(f))$ in the case of equality take $f = x^nf_n$ and $D = x^m u\frac{\partial}{\partial x_i}$ where $pdeg_x(f_n) = 0$ and $u \in L(x)$. Since $D(f) = x^n D(f_n) + nx^{n-1}f_n D(x)$ only $x^n D(f_n)$ can contain $x^{n+m}$ and we should show that $l_x(D(f_n)) = x^m D_m(f_n)$ where $D_m = u\frac{\partial}{\partial x_i}$. It can be done exactly as above by reduction first to the case when $pdeg(f_n) = 1$ and then to the case when $deg(f_n) = 1$. $\Box$

\begin{lm}\label{l1}
Let $D$ be a derivation of $P$ and $x$ be the minimal element of $S(D)$.
If $D$ is locally nilpotent then so is $l_x(D)$.
\end{lm}
\Proof If $l_x(D)$ is not locally nilpotent then there exists
$x_i$ such that  $l_x(D)^k(x_i)\neq 0$ for all $k\geq 0$. Put
$a=l_x(D)(x_i)$. Note that $x\notin S(a)$ and $l_x(a)=a$. Using
this and Proposition \ref{pr2}, we get \bes
l_x(D^k(a))=l_x(D)(l_x(D^{k-1}(a)))=\ldots=l_x(D)^{k-1}(l_x(D)(a))=l_x(D)^k(a)\neq
0. \ees Consequently, $D$ is not locally nilpotent. $\Box$

\begin{pr}\label{pr3}
Let $D$ be a derivation of $P$ of the form \bes
D=D_0+xD_1+\ldots+x^{m-1}D_{m-1}+x^m\frac{\partial}{\partial x_1},
\ \ x\notin S(D_i), \ \ 0\leq i\leq m-1, \ees where $x$ is the
minimal element of $S(D)$. Let $f$ be an element of $P$ such that
$x_1\notin S(f)$. Then \bes pdeg_xD(f)\leq m-1+pdeg_xf. \ees This
inequality becomes an equality iff $D'(l_x(f))\neq 0$, where
$D'=D_{m-1}+m x\frac{\partial}{\partial x_1}$, and in this case
\bes l_x(D(f))=D'(l_x(f)). \ees
\end{pr}
\Proof The same considerations as in the proof of Proposition
\ref{pr2} show that \bes pdeg_x(x^m\frac{\partial}{\partial
x_1}(f))\leq m-1+pdeg_xf \ees and if $\frac{\partial}{\partial
x_1}(l_x(f))\neq 0$ then \bes l_x(x^m\frac{\partial}{\partial
x_1}(f))=m l_x(x\frac{\partial}{\partial x_1}(f)). \ees

Note that $D=D^*+x^m\frac{\partial}{\partial x_1}$ and
$pdeg_x(D^*)\leq m-1$. So applying Proposition \ref{pr2}, we can
complete the proof of Proposition \ref{pr3}.  $\Box$

\begin{lm}\label{l2}
Let $D$ be a locally nilpotent derivation of $P$ of the form
\bes
D=D_0+xD_1+\ldots+x^{m-1}D_{m-1}+x^m\frac{\partial}{\partial x_1}, \ \ x\notin S(D_i), \ \ 0\leq i\leq m-1,
\ees
where $x$ is the minimal element of $S(D)$. If $x\neq x_1$ then $D_{m-1}+m x\frac{\partial}{\partial x_1}$
is also locally nilpotent.
\end{lm}
\Proof Assume that $D'=D_{m-1}+m x\frac{\partial}{\partial x_1}$
is not locally nilpotent. Then there exists $x_i$ such that
$D'^k(x_i)\neq 0$ for all $k\geq 0$. We put $a=D'^2(x_i)$. It is
not difficult to show that $x_1, x \notin S(a)$. So $l_x(a)=a$.
Using this and Proposition \ref{pr3}, we get \bes
l_x(D^k(a))=D'(l_x(D^{k-1}(a)))=\ldots=D'^k(a)\neq 0. \ees
Consequently, $D$ is not locally nilpotent. $\Box$

\begin{lm}\label{l3}
Let $D$ be a multihomogeneous derivation of $P=P\langle x_1,x_2\rangle$ and $mdeg(D)= (m_1,m_2)$.
If $m_i\geq 0$ for $i=1,2$ then $D$ is not locally nilpotent.
\end{lm}
\Proof Let $D$ be a counterexample to the lemma with the minimal
$\deg(D)$. By Proposition \ref{pr1}, we can also assume  that $D$
is $p$-homogeneous. Let $x$ be the minimal element of $S(D)$. By
Lemma \ref{l1}, it follows that $l_x(D)$ is also locally
nilpotent. Put $mdeg(l_x(D))=(n_1,n_2)$. We can assume that
$n_1=-1$ since $\deg(l_x(D))<\deg(D)$. Then $l_x(D) = \a x_2^{n_2}
\frac{\partial}{\partial x_1}$.

If $x=x_1$ then $D$ contains a summand $l_x(D)= \a x_1^{m_1+1}
x_2^r \frac{\partial}{\partial x_1}$. In this case $D$ induces a
nonzero locally nilpotent derivation of the polynomial algebra $k[x_1,x_2]$ with the
same multidegree. It is impossible (see, for example \cite{Essen},
p. 91).

So $x \neq x_1$. If $x = x_2$ then $m_1 = -1$. So $x>x_2$ and $D$
can be written as in Lemma \ref{l2}. By Lemma \ref{l2}, it follows
that $D'=D_{m-1}+m x\frac{\partial}{\partial x_1}$ is a nonzero
locally nilpotent derivation. Note that $pdeg(D')=0$ and $D'$ is
$p$-homogeneous. Therefore $D'$ is a derivation of the free Lie
algebra $g$  generated by  $x_1,x_2$. Obviously, $exp(D')$ gives a
nonlinear automorphism of $g$. But all automorphisms of $g$ are
linear \cite{Cohn2}. $\Box$

\section{The main results}

\hspace*{\parindent}

Recall that a derivation of the free Poisson algebra
$P\langle x_1,x_2,\ldots,x_n\rangle$ of the form (\ref{f3}) is called
triangular if $f_i\in P\langle x_{i+1},x_{i+2},\ldots,x_n\rangle$ for any $i$.
It is clear that every triangular derivation is locally nilpotent.
A derivation $D$ of $P\langle x_1,x_2,\ldots,x_n\rangle$ is called triangulable
if there exists an automorphism $\varphi$ such that
$\varphi^{-1}D\varphi$ is triangular. R.\, Rentschler proved
\cite{Rentschler}  that the locally nilpotent derivations of
polynomial algebras in two variables over a field of
characteristic $0$ are triangulable. H.\,Bass gave \cite{Bass2} an
example of a nontriangulable derivation of polynomial algebras in
three variables.

\begin{th}\label{t1}
Let $D$ be a locally nilpotent derivation of $P=P\langle x_1,x_2\rangle$. Then
there exist a tame automorphism $\varphi$ of $P$ and $f(x_2)\in
k[x_2]$ such that
$\varphi^{-1}D\varphi=f(x_2)\frac{\partial}{\partial x_1}$.
\end{th}
\Proof Denote by $I$ the ideal of $P$ generated by $\{x_1,x_2\}$.
Then $P/I\cong k[x_1,x_2]$ and $D$ induces a locally nilpotent
derivation $D'$ of $k[x_1,x_2]$. By Rentschler's theorem
\cite{Rentschler}, there exists a tame automorphism $\psi$ of
$k[x_1,x_2]$ and $f(x_2)\in k[x_2]$ such that
$\psi^{-1}D'\psi=f(x_2)\frac{\partial}{\partial x_1}$. 
Denote by $\varphi$ the extension of $\psi$ to $P$
such that $\varphi|_{k[x_1,x_2]}=\psi$. Replacing $D$ by
$\varphi^{-1}D\varphi$ we can assume that
$D'=f(x_2)\frac{\partial}{\partial x_1}$. Then \bes
D=(f(x_2)+a)\frac{\partial}{\partial
x_1}+b\frac{\partial}{\partial x_2}, \ees where $a, b \in I$.

We would like to show that $a = b = 0$. Assume it is not the case.
Consider $\deg_{x_1}$ and the corresponding highest homogeneous
derivation $R$ which is locally nilpotent by Proposition
\ref{pr1}. But $R = c\frac{\partial}{\partial x_1} +
d\frac{\partial}{\partial x_2}$ where $c, d \in I$ and either $c$
or $d$ is not zero. So $R$ cannot be locally nilpotent by Lemma 3.
$\Box$

\begin{co}\label{c1}
Let $D$ be a locally nilpotent derivation of $P=P\langle x_1,x_2\rangle$.
Then $D\{x_1,x_2\}=0$.
\end{co}
\Proof If $D$ is triangular then $D\{x_1,x_2\}=0$. Note that
$\varphi\{x_1,x_2\}=\a \{x_1,x_2\}$ for every tame automorphism
since it is true for every elementary automorphism.  $\Box$

\begin{th}\label{t2}
Automorphisms of free Poisson algebras in two variables over a
field of characteristic zero are tame.
\end{th}
\Proof Let $\theta$ be an arbitrary automorphism of
$P=P\langle x_1,x_2\rangle$. Then $\theta$ induces an automorphism $\psi$ of
$k[x_1,x_2]$. Denote by $\varphi$ the extension of $\psi$ to $P$
such that $\varphi|_{k[x_1,x_2]}=\psi$. By Jung's theorem
\cite{Jung}, $\psi$ and $\varphi$ are tame. Changing  $\theta$ to
$\theta\varphi^{-1}$ we can assume that $\theta$ induces the
identical automorphism of $k[x_1,x_2]$. Then, \bes
\theta(x_1)=x_1+a, \ \ \theta(x_2)=x_2+b, \ \ a, b\in I, \ees
where $I$ is the ideal of $P$ generated by $\{x_1,x_2\}$.

For every $h\in k[x]$ denote by $D_h$ a derivation of $P$ defined
by $D_h(x_1+a)= h(x_2+b)$, $D_h(x_2+b)=0$. This derivation is
locally nilpotent.

Now,  \bes D_h = (h(x_2) +(h(x_2+b) - h(x_2)) -
D_h(a))\frac{\partial}{\partial x_1}-
D_h(b)\frac{\partial}{\partial x_2} \ees since $D_h(x_1) = h(x_2)
+(h(x_2+b) - h(x_2)) - D_h(a)$ and $D_h(x_2) = - D_h(b)$.

The ideal $I$ is invariant under every derivation. Hence $h(x_2+b)
- h(x_2)) - D_h(a), D(b) \in I$. Since $D_h$ is locally nilpotent
it is possible only if $D_h(b) = h(x_2+b) - h(x_2)) - D_h(a) = 0$
(see the proof of Theorem 1). Therefore $D_h(x_1)=h(x_2),
 \ \ D_h(x_2)=0$ and $D_h(a)=h(x_2+b)-h(x_2)$.

Put $h = x$. Then $D_x(a)=b$. Note that $\deg\,D_x(a)\leq
\deg\,a$ since $D_x(x_1) = x_2$ and $D_x(x_2) = 0$. 
So $\deg\,b\leq \deg\,a$. We can exchange $x_1$
and $x_2$ in the definition of $D_h$, so $\deg\,a\leq
\deg\,b$ and $\deg\,a=\deg\,b$. Of course, $\deg\,a=\deg\,b\geq
2$ since $a, b \in I$.

We now put $h =x^2$. Then $D_h(a)= 2 x_2 b+ b^2$. Note that in
this case $\deg\,D_h(a)\leq \deg\,a+1$ since $D_h(x_1) = x_2^2$ 
and $D_h(x_2) = 0$. Consequently,
$\deg\,a+1\geq 2 \deg\,b=2 \deg\,a$, and $\deg\,a\leq 1$. This
contradiction gives $a=0$ and $b=0$. $\Box$

\begin{co}\label{c2}
Let $\varphi$ be an arbitrary automorphism of $P=P\langle x_1,x_2\rangle$.
Then $\varphi\{x_1,x_2\}=\a \{x_1,x_2\}$, where $0\neq \a \in k$.
\end{co}

So every automorphism of $P\langle x_1,x_2\rangle$ preserves $\{x_1,x_2\}$ up to the 
proportionality. An analogue
of this result for free associative algebras is also true, i.e., every automorphism
of the free associative algebra $k<x_1,x_2>$ in the variables $x_1,x_2$
preserves the commutator $[x_1,x_2]$ up to the
proportionality. Moreover, the so called commutator test theorem \cite{Dicks}
says that any endomorphism of $k<x_1,x_2>$ which preserves $[x_1,x_2]$ is an automorphism.

\begin{prob}\label{prob1}
Is any endomorphism of the free Poisson algebra $P\langle x_1,x_2\rangle$ over
a field of characteristic $0$ which preserves $\{x_1,x_2\}$ an automorphism?
\end{prob}

Note that the positive answer to Problem \ref{prob1} implies the Jacobian Conjecture 
for $k[x_1,x_2]$ \cite{Essen}.\\

It is well known \cite{Czer,Makar} that $Aut\,k[x_1,x_2]\cong Aut\,k<x_1,x_2>$,
where  $k<x_1,x_2>$ is the free associative algebra generated by  $x_1,x_2$.

\begin{co}\label{c3} Let $k$ be a field of characteristic zero. Then,
\bes
Aut\,k[x_1,x_2]\cong Aut\,k<x_1,x_2>\cong Aut\,P\langle x_1,x_2\rangle.
\ees
\end{co}

This isomorphism is also interesting in the context of paper \cite{Belov-Kontsevich}
since $k<x_1,x_2>$ is a deformation quantization of $P\langle x_1,x_2\rangle$ and because it shows that
the group $Aut\,P\langle x_1,x_2\rangle$ has a nice representation as a free amalgamated product of its subgroups
(see, for example \cite{Cohn}).

\bigskip

\begin{center}
{\bf\large Acknowledgments}
\end{center}

\hspace*{\parindent}

The authors wish to thank several institutions which supported them while they were 
working on this project: 
Max-Planck Institute
f\"ur Mathematik (the first and the third authors), Department of Mathematics of Wayne State University in Detroit (the third author), and Instituto de Matem\'{a}tica e Estat\'{i}stica da Universidade de S\~{a}o Paulo (the first author).


\begin{thebibliography}{99}

\bibitem{Bass2} H.\,Bass, A non-triangular action of $G_a$ on $A^3$, J. of Pure and
Appl. Algebra, 33(1984), no.1, 1--5.

\bibitem{Belov-Kontsevich} A.\,Belov-Kanel, M.\,Kontsevich, Automorphisms of the Weyl Algebra,
Letters in Mathematical Physics, 74 (2005), 181--199.


\bibitem{Berg} G.\,M.\,Bergman, Centralizers in free associative algebras,
Trans. Amer. Math. Soc., 137 (1969), 327--344.

\bibitem{Cohn2} P.\,M.\,Cohn, Subalgebras of free associative algebras,
Proc. London Math. Soc., 56 (1964), 618--632.

\bibitem{Cohn} P.\,M.\,Cohn, Free rings and their relations, 2nd Ed.,
Academic Press, London, 1985.

\bibitem{Czer} A.\,G.\,Czerniakiewicz, Automorphisms of a free associative
algebra of rank 2, I, II, Trans. Amer. Math. Soc., 160 (1971),
393--401; 171 (1972), 309--315.

\bibitem{Dicks} W.\,Dicks, A commutator test for two elements to generate the free algebra of rank two,
Bull. London Math. Soc., 14 (1982), 48--51.

\bibitem{Essen} A.\,van den Essen, Polynomial automorphisms and the
Jacobian conjecture, Progress in Mathematics, 190, Birkhauser verlag,
Basel, 2000.


\bibitem{Jung} H.\,W.\,E.\,Jung, Uber ganze birationale Transformationen
der Ebene, J. reine angew. Math., 184 (1942), 161--174.


\bibitem{Kulk} W.\,van der Kulk, On polynomial rings in two variables,
Nieuw Archief voor Wiskunde, (3)1 (1953), 33--41.


\bibitem{Makar} L.\,Makar-Limanov, The automorphisms of the free algebra with two generators, Funksional.
Anal. i Prilozhen. 4(1970), no.3, 107-108; English translation: in
Functional Anal. Appl. 4 (1970), 262--263.

\bibitem{MakarU2}
L.~Makar-Limanov, U.~U.~Umirbaev, Centralizers in free Poisson algebras,
  Proc. Amer. Math. Soc. 135 (2007), no. 7,  1969--1975.

\bibitem{Nagata} M.\,Nagata, On the automorphism group of $k[x,y]$,
Lect. in Math., Kyoto Univ., Kinokuniya, Tokio, 1972.

\bibitem{Rentschler} R.\,Rentschler, Operations du groupe additif sur le plan, C.\,R. Acad. Sci. Paris, 267 (1968), 384--387.

\bibitem{Shest} I.\,P.\,Shestakov, Quantization of Poisson superalgebras
and speciality of Jordan Poisson superalgebras, Algebra i logika,
32(1993), no. 5, 571--584; English translation: in Algebra and
Logic, 32(1993), no. 5, 309--317.

\bibitem{Umi24}
I.~P.~Shestakov and U.~U.~Umirbaev, Poisson brackets and two
generated subalgebras of rings of polynomials, Journal of the
American Mathematical Society, 17 (2004), 181--196.

\bibitem{Umi25}
I.~P.~Shestakov and U.~U.~Umirbaev, Tame and wild automorphisms of
rings of polynomials in three variables, Journal of the American
Mathematical Society, 17 (2004), 197--227.

\bibitem{Um19}
 U.~U.~Umirbaev,  The Anick automorphism of free associative algebras,
 J. Reine Angew. Math. 605 (2007), 165--178.

\end{thebibliography}
\end{document}